\documentclass[12pt,a4paper]{amsart}

\title[Abelian subgroup structure]{Abelian subgroup structure of square complex groups 
and arithmetic of quaternions}
\author{Diego Rattaggi}
\address{Universit\'e de Gen\`eve,
Section de math\'ematiques,
2--4 rue du Li\`evre, CP 64, 
CH--1211 Gen\`eve 4, Switzerland.}
\email{rattaggi@math.unige.ch}
\author{Guyan Robertson}
\address{School of Mathematics and Statistics, University of Newcastle, NE1 7RU, U.K.}
\email{a.g.robertson@newcastle.ac.uk}
\date{\today}

\usepackage{amsmath,amssymb}
\usepackage{amscd}

\input{prepictex}
\input{pictex}
\input{postpictex}
\chardef\bslash=`\\ 

\makeatletter
\def\verbatim{\interlinepenalty\@M \@verbatim
  \leftskip\@totalleftmargin\advance\leftskip2pc
  \frenchspacing\@vobeyspaces \@xverbatim}
\makeatother
\newtheorem{theorem}{Theorem}[section]
\newtheorem{corollary}[theorem]{Corollary}
\newtheorem{lemma}[theorem]{Lemma}
\newtheorem{proposition}[theorem]{Proposition}

\theoremstyle{definition}
\newtheorem{example}[theorem]{Example}
\newtheorem{definition}[theorem]{Definition}

\newtheorem{notation}[theorem]{Notation}
\newcommand{\ovl}{\overline}

\newcommand{\cl}[1]{{\mathcal{#1}}}
\newcommand{\bb}[1]{{\mathbb{#1}}}

\newcounter{picture}

\newcommand{\SO}{{\text{\rm{SO}}}}

\begin{document}
\def\Proof  {{\bf Proof.}\par}
\def\bs {\backslash}
\def\trace   {{\text{\rm trace}}}
\maketitle

\section{Introduction}

A square complex is a 2-complex formed by gluing squares together.
This article is concerned with the fundamental group $\Gamma$
of certain square complexes of nonpositive curvature, related to 
quaternion algebras. The abelian subgroup structure of $\Gamma$
is studied in some detail. Before outlining the results, it is
necessary to describe the construction of $\Gamma$.

In \cite[Section 3]{moz}, there is constructed a lattice subgroup 
$\Gamma = \Gamma_{p,l}$ of $G = PGL_2(\bb Q_p) \times
PGL_2(\bb Q_l)$, where $p,l \equiv 1 \pmod 4$ are two distinct
primes. This restriction was made because $-1$ has a square root in $\bb Q_p$
if and only if $p \equiv 1 \pmod 4$,  
but the construction of $\Gamma$ is generalized in \cite[Chapter~3]{rat}
to all pairs $(p,l)$ of distinct odd primes.

The affine building $\Delta$ of $G$ is a product
of two homogeneous trees of degrees $(p+1)$ and $(l+1)$ respectively.
The group $\Gamma$ is a finitely presented torsion free group which acts freely and
transitively on the vertices of $\Delta$,
with a finite square complex as quotient $\Delta / \Gamma$.

Here is how $\Gamma$ is constructed. Let
$$\bb H(\bb Z)=\{ x=x_0+x_1i+x_2j+x_3k ; x_0, x_1, x_2, x_3 \in \bb Z\}$$
be the ring of integer quaternions where $i^2 = j^2 = k^2 = -1$, $ij = -ji = k$.
Let $\ovl x = x_0 - x_1 i - x_2 j - x_3 k$ be the conjugate of $x$,
and $|x|^2 = x \ovl x = x_0^2 + x_1^2 + x_2^2 + x_3^2$ its norm.

Let $c_p, d_p \in \bb Q_p$ and $c_l, d_l \in \bb Q_l$ be elements such that
$c_p^2 + d_p^2 + 1 = 0$, $c_l^2 + d_l^2 + 1 = 0$.
Such elements exist by Hensel's Lemma and \cite[Proposition 2.5.3]{dsv}.
We can take
$d_p = 0$, if $p \equiv 1 \pmod 4$, and
$d_l = 0$, if $l \equiv 1 \pmod 4$.
Define
$$\psi : \bb H(\bb Z) - \{ 0 \} \to PGL_2(\bb Q_p) \times PGL_2(\bb Q_l)$$
 by
\begin{equation}\label{psimap}
\begin{split}
\psi(x)=\bigg(
&\begin{pmatrix}
x_0 + x_1 c_p + x_3 d_p & -x_1 d_p + x_2 + x_3 c_p \\
-x_1 d_p - x_2 + x_3 c_p  & x_0 - x_1 c_p - x_3 d_p \\
\end{pmatrix}, \\
&\begin{pmatrix}
x_0 + x_1 c_l + x_3 d_l & -x_1 d_l + x_2 + x_3 c_l \\
-x_1 d_l - x_2 + x_3 c_l  & x_0 - x_1 c_l - x_3 d_l \\
\end{pmatrix}
\bigg).
\end{split}
\end{equation}
This formula abuses notation by identifying an element of 
$PGL_2(\bb Q_p) \times PGL_2(\bb Q_l)$ with its representative in
$GL_2(\bb Q_p) \times GL_2(\bb Q_l)$.

Note that $\psi(xy) = \psi(x) \psi(y)$, $\psi(\lambda x) = \psi(x)$, if $\lambda \in \bb Z - \{ 0 \}$,
and $\psi(x)^{-1} = \psi(\ovl x )$.
Moreover the inverse image under $\psi$ of the identity element in $PGL_2(\bb Q_p) \times PGL_2(\bb Q_l)$ is precisely 
$$\bb Z - \{ 0 \} = \{ x \in \bb H(\bb Z); x_0 \ne 0, x_1 = x_2 = x_3 = 0 \}\,.$$
Let 
\begin{align*}
\tilde\Gamma=\{x \in \bb H(\bb Z) \,; 
\quad& |x|^2=p^rl^s, r,s \ge 0\,;  \notag \\
&x_0 \text{ odd}, x_1, x_2, x_3 \text{ even}, \text{ if } |x|^2 \equiv 1 \pmod 4\,; \notag \\
&x_1 \text{ even}, x_0, x_2, x_3 \text{ odd}, \text{ if } |x|^2 \equiv 3 \pmod 4
\}.
\end{align*}
Then $\Gamma=\psi(\tilde\Gamma$) is a torsion free cocompact lattice in $G$.
Let
\begin{equation*}
\tilde A =\{x \in \tilde\Gamma ; x_0>0, |x|^2=p\}, \quad
\tilde B =\{y \in \tilde\Gamma ; y_0>0, |y|^2=l\}\,.
\end{equation*}
Then $\tilde A$ contains $p+1$ elements and $\tilde B$ contains $l+1$ elements,
by a result of Jacobi \cite [Theorem 2.1.8]{lub}.
The images $A = \psi(\tilde A), B = \psi(\tilde B)$ of $\tilde A, \tilde B$ in $\Gamma$ generate
free groups $\Gamma_p = \langle A \rangle = \langle a_1, \ldots, a_{\frac{p+1}{2}} \rangle$, 
$\Gamma_l = \langle B \rangle = \langle b_1, \ldots, b_{\frac{l+1}{2}} \rangle$ 
of ranks $(p+1)/2$, $(l+1)/2$ respectively and $\Gamma$
itself is generated by $A \cup B$.
The 1-skeleton of $\Delta$ is the Cayley graph of $\Gamma$ relative to this set
of generators.
The group $\Gamma$ has a finite presentation with generators
$\{ a_1, \ldots, a_{\frac{p+1}{2}} \} \cup \{ b_1, \ldots, b_{\frac{l+1}{2}} \}$
and $(p+1)(l+1)/4$ relations of the form $ab = \tilde b \tilde a$, where $a, \tilde a \in A$,
$b, \tilde b \in B$.
In fact, given any $a \in A$, $b \in B$, there are unique elements $\tilde a \in A$, 
$\tilde b \in B$ such that $ab = \tilde b \tilde a$.
This follows from a special case of Dickson's factorization property for integer quaternions
(\cite[Theorem~8]{di}).
\begin{proposition} \label{dickson} (\cite{di})
Let $x \in \tilde \Gamma$ such that $|x|^2 = pl$.
Then there are uniquely determined $z, \tilde z \in \tilde A$, $y, \tilde y \in \tilde B$
such that $zy, \tilde y \tilde z = \pm x$.
\end{proposition}

It is worth noting that $zy\not= \tilde y \tilde z$ in general, as demonstrated by the following example. 

\begin{example}
Let $p=3$, $l=5$ and $x = 1+2i+j+3k$. Then
$(1-j+k)(1+2i) = x$ and $(1-2k)(1-j-k) = -x$. 
\end{example}

We can now outline the contents of this article. A fundamental fact, upon which much else depends,
is that $\Gamma$ is {\em commutative transitive}, in the sense that the relation of commutativity is transitive on non-trivial elements of $\Gamma$. In particular $\Gamma$ cannot contain
a subgroup isomorphic to $F_2 \times F_2$, where $F_2$ denotes the free group of rank $2$.
Furthermore, $\Gamma$ is a {\it CSA-group},
i.e.\ all its maximal abelian subgroups $\Gamma_0$ satisfy $g \Gamma_0 g^{-1} \cap \Gamma_0 = \{1\}$ for all $g\in \Gamma -\Gamma_0$.

Every nontrivial element $\gamma\in \Gamma$ is the image under $\psi$ of a quaternion of the form
$x_0 + z_0 (c_1 i + c_2 j + c_3 k)\,$
where $c_1, c_2, c_3\in\bb Z$ are relatively prime.
The element $\gamma$ is contained in a unique maximal abelian subgroup $\Gamma_0$ and
the integer $n=n(\Gamma_0)= c_1^2+c_2^2+c_3^2$ depends only on $\Gamma_0$ rather than the particular
choice of $\gamma$.
We define a class of maximal abelian subgroups of $\Gamma$ isomorphic to $\bb Z^2$,
which we call period subgroups,
and which are characterized by the condition $\left(\frac{-n}{p}\right)=\left(\frac{-n}{l}\right)=1$.
Every maximal abelian subgroup $\Gamma_0 \cong \bb Z^2$ is conjugate in $\Gamma$ to a period subgroup 
and, as the name suggests, period subgroups are closely related to periodic tilings of the plane.
On the other hand, some maximal abelian subgroups of $\Gamma$ are isomorphic to $\bb Z$, and we show how to construct these. Several explicit examples and counterexamples are included.

\section{The CSA property}

Let $\tau : \bb H(\bb Q) - \bb Q \to \bb P^2(\bb Q)$ be defined by $\tau(x) = \bb Q(x_1, x_2, x_3)$, which is 
a line in $\bb Q^3$ through $(0,0,0)$. 
By \cite[Section 3]{moz} , two quaternions $x, y \in \bb H(\bb Q) - \bb Q$
commute if and only if $\tau(x)=\tau(y)$.
This directly implies the following lemma, 
which in turn has
Proposition~\ref{commtrans} as a consequence,
see also \cite[Chapter~3]{rat}.

\begin{lemma}\label{comm}
Elements $x, y \in \tilde \Gamma$ commute if and only if their images $\psi(x), \psi(y) \in \Gamma$ commute.
\end{lemma}

A group is said to be {\em commutative transitive} if the relation of commutativity is transitive on its non-trivial elements.

\begin{proposition}\label{commtrans}
The group $\Gamma$ is commutative transitive.
\end{proposition}

Wise has asked in \cite[Problem~10.9]{wis} whether the fundamental group of any 
nonelementary complete square complex contains a subgroup isomorphic to $F_2 \times F_2$.
We can give a negative answer of this question, since our group $\Gamma$ belongs to this class of fundamental groups,
and it is a direct consequence of Proposition~\ref{commtrans} that $\Gamma$ does not contain a $F_2 \times F_2$ subgroup. 
In fact, since $\Gamma$ is torsion free, 
and a (free) abelian subgroup of $\Gamma$ has rank $\le 2$ \cite[Lemma 3.2]{pra}, we have a more precise result.

\begin{corollary}\label{corwise}
The only nontrivial direct product subgroup of $\Gamma$ is
$\bb Z \times \bb Z = \bb Z^2\,.$
\end{corollary}

If $\gamma = \psi(x) \in \Gamma - \{ 1 \}$ then the centralizer $\Gamma_0 = Z_{\Gamma}(\gamma)$ is the unique maximal abelian subgroup 
of $\Gamma$ containing $\gamma$. Moreover $\Gamma_0$ is determined by $\tau(x)$, independent of the choice of $x$.

As described in \cite[Remark~4]{mr},
a group is commutative transitive if and only if the centralizer of any non-trivial element is abelian.
A third equivalent condition (called {\it SA-property} in \cite{mr})
is proved for $\Gamma$ in the following lemma.
It is used to show in Proposition~\ref{csa} that $\Gamma$ is a {\it CSA-group},
i.e.\ all its maximal abelian subgroups are malnormal,
where a subgroup $\Gamma_0$ of $\Gamma$ is {\it malnormal} (or {\it conjugate separated}) if 
$g \Gamma_0 g^{-1} \cap \Gamma_0 = \{1\}$ for all $g\in \Gamma -\Gamma_0$.
Any CSA-group is commutative transitive, but the converse is not true, see \cite{mr}.

\begin{lemma}\label{maxabelian}
If $\Gamma_1$ and $\Gamma_2$ are maximal abelian subgroups of $\Gamma$
and $\Gamma_1\not=\Gamma_2$ then $\Gamma_1 \cap \Gamma_2 = \{1\}$.
\end{lemma}

\begin{proof}
Suppose that there exists a nontrivial element $\gamma \in \Gamma_1 \cap \Gamma_2$.
If $\gamma_i\in \Gamma_i - \{1\}$, $i=1, 2$, then 
$\gamma \gamma_1 = \gamma_1 \gamma$ and $\gamma \gamma_2 = \gamma_2 \gamma$ 
which implies $\gamma_1 \gamma_2 = \gamma_2 \gamma_1$ by Proposition~\ref{commtrans}.
Since $\Gamma_1$, $\Gamma_2$ are maximal abelian, $\Gamma_1=\Gamma_2$.
\end{proof}
 
It is well known
that there is a (surjective) homomorphism
\[
\theta : \bb H(\bb Q) - \{0\} \to \SO_3(\bb Q)
\]
defined by $\theta(y) x = yxy^{-1}$, for $x=x_1 i + x_2 j + x_3 k \in \bb H(\bb Q)$
identified with $(x_1, x_2, x_3) \in \bb Q^3$.

If $y\in \bb H(\bb Q) - \bb Q$ then the axis of rotation of $\theta(y)$ is $\tau(y)$.
This is an immediate consequence of the fact that 
$$\theta(y)(y-y_0)=y(y-y_0)y^{-1} = y-y_0\,.$$
Moreover the angle of rotation is $2\alpha$ where $\cos \alpha = \frac{y_0}{|y|}$
\cite[Chapitre I, \S 3]{v}. In particular, the angle of rotation is a multiple of $\pi$
only if $y_0=0$.

\begin{lemma}\label{conj}
(a) Suppose that $x, y \in \bb H(\bb Q) - \bb Q$ and $y_0 \not= 0$. 
Then  $yxy^{-1}$ commutes with $x$ if and only if $y$ commutes with $x$.

(b) If $a, b \in \Gamma$, then  $bab^{-1}$ commutes with $a$
if and only if $b$ commutes with $a$.
\end{lemma}

\begin{proof} (a)
If $yxy^{-1}$ commutes with $x$, then the rotations $\theta(yxy^{-1})$ and $\theta(x)$ have the same axis.  
However, the axis of $\theta(yxy^{-1}) = \theta(y)\theta(x)\theta(y)^{-1}$ is $\theta(y)\tau(x)$.
Therefore $\theta(y) \tau(x)=\tau(x)$ : in other words $\theta(y) (x_1,x_2,x_3)=\pm(x_1,x_2,x_3)$.
Now if $\theta(y) (x_1,x_2,x_3)=-(x_1,x_2,x_3)$ then $\theta(y)$ is a rotation of angle $\pi$,
with axis perpendicular to $(x_1,x_2,x_3)$. This cannot happen since $y_0 \not= 0$. 
Therefore $\theta(y)$ has axis $\tau(x)$. That is $\tau(y)=\tau(x)$, and consequently
$y$ commutes with $x$. The converse is clear.

(b) If $a = 1$ or $b = 1$, the statement is obvious. 
If $a, b \in \Gamma - \{1\}$ and $bab^{-1}$ commutes with $a$, then representatives $x,y$ for
$a, b$ in $\bb H(\bb Q) - \bb Q$ have nonzero real parts and satisfy the same relation, by Lemma~\ref{comm}. 
The assertion follows from (a). Again, the converse is clear.
\end{proof}

\begin{proposition}\label{csa}
$\Gamma$ is CSA.
\end{proposition}

\begin{proof}
Suppose that $\Gamma_0$ is a maximal abelian subgroup of $\Gamma$ and that
$b \in \Gamma$, with $b \Gamma_0 b^{-1} \cap \Gamma_0 \not= \{1\}$.
We must show that $b \in \Gamma_0$.

By Lemma~\ref{maxabelian}, $b \Gamma_0 b^{-1} = \Gamma_0$.  Let $a\in \Gamma_0$.
Then $b a b^{-1}$ commutes with $a$ and so, by Lemma~\ref{conj}, $b$ commutes with $a$.
Since $\Gamma_0$ is maximal abelian, $b \in \Gamma_0$. 
\end{proof}

We now recall the following known result.

\begin{lemma}\label{mrcsa}
(a) (\cite[Proposition~9(5)]{mr}) A non-abelian CSA-group has no
non-abelian solvable subgroups.

(b) (\cite[Proposition~10(3)]{mr}) Subgroups of CSA-groups are CSA.
\end{lemma}

\begin{corollary}\label{titsalternative}
Let $a \in \Gamma_p - \{ 1 \}$ and $b \in \Gamma_l - \{ 1 \}$.
Then either $\langle a, b \rangle \cong \bb Z^2$ or $\langle a, b \rangle$
contains a free subgroup of rank $2$.
\end{corollary}

\begin{proof}
If $a, b$ commute, then $\langle a, b \rangle \cong \bb Z^2$, since $\Gamma$ is torsion free
and $\langle a, b \rangle$ is not cyclic.
Assume that $a, b$ do not commute.
We will show that $\langle a, b \rangle$ is not virtually solvable.
The Tits Alternative for finitely generated linear groups (see \cite{ti})
then implies that $\langle a, b \rangle$ contains a free subgroup of rank $2$.
Note that $\Gamma$ is linear, see \cite[Section~3.2]{rat}
for an explicit injective homomorphism $\Gamma \to SO_3(\bb Q)$.
Let $U$ be a finite index subgroup of $\langle a, b \rangle$,
in particular there are $r,s \in \bb N$ such that $a^r, b^s \in U$.
The elements $a^r$ and $b^s$ do not commute since otherwise also $a$ and $b$ would commute by
Proposition~\ref{commtrans}. It follows that $U$ is not abelian.
By Proposition~\ref{csa} and Lemma~\ref{mrcsa}(b), $\langle a, b \rangle$ is CSA.
Lemma~\ref{mrcsa}(a) shows that $U$ is not solvable.
\end{proof}

\bigskip

\section{Maximal abelian subgroups and period subgroups.}\label{persg}

\bigskip
Recall that the group $\Gamma$ acts freely and transitively on the vertex set of the
affine building $\Delta$ of $PGL_2(\bb Q_p) \times PGL_2(\bb Q_l)$. The building $\Delta$
is a product of two homogeneous trees and the apartments (maximal flats) in $\Delta$ are
copies of the Euclidean plane tesselated by squares. 

\bigskip

\begin{notation}
If $n$ is an integer and $p$ is an odd prime, then the {\bf Legendre symbol} is
\begin{equation*} \left(\frac{n}{p}\right)=
\begin{cases}
0 & \text{if $p \mid n$},\\
1 & \text{if $p \nmid n$ and $n$ is a square mod $p$},\\
-1 & \text{if $p \nmid n$ and $n$ is not a square mod $p$}.  
\end{cases}
\end{equation*}
\end{notation}

Any element of $\Gamma-\{1\}$ is the image under $\psi$ of a quaternion of the form
\begin{equation}\label{x}
x = x_0 + z_0 (c_1 i + c_2 j + c_3 k)\, ,
\end{equation}
where $c_1, c_2, c_3\in\bb Z$ are relatively prime, $z_0 \ne 0$, $(c_1, c_2, c_3) \ne (0,0,0)$, and
$$|x|^2 = x_0^2 + (c_1^2 + c_2^2 + c_3^2) z_0^2 = p^r l^s, r,s \ge 0\,.$$

Recall that $\tau(x)=\bb Q(c_1, c_2, c_3)\in \bb P^2(\bb Q)$ and 
recall that elements $\psi(x), \psi(y) \in \Gamma-\{1\}$ commute if and
only if $\tau(x)=\tau(y)$.
Moreover the centralizer $\Gamma_0=Z_{\Gamma}(\psi(x))$ is the unique maximal abelian subgroup of
$\Gamma$ containing $\psi(x)$. Let
$$n(x) = n(\psi(x)) = n(\Gamma_0)= c_1^2+c_2^2+c_3^2\,.$$

An abelian subgroup of $\Gamma$ has rank $\le 2$ \cite[Lemma 3.2]{pra}.
Since $\Gamma$ is torsion free,
a nontrivial abelian subgroup $\Gamma_0$ of $\Gamma$ is isomorphic to either
$\bb Z$ or $\bb Z^2$. 
If $\Gamma_0 \cong \bb Z^2$ then there is a unique apartment $\cl A_{\Gamma_0}$ which is stabilized by
$\Gamma_0$ \cite[6.8]{pra}, and $\Gamma_0$ acts cocompactly by translation on
this apartment. We call $\cl A_{\Gamma_0}$ a {\em periodic} apartment.

\begin{definition}
A maximal abelian subgroup $\Gamma_0 \cong \bb Z^2$ will be called a {\bf period subgroup}
if the apartment $\cl A_{\Gamma_0}$ contains 
the vertex $O$ of $\Delta$ whose stabilizer in $G$ is $PGL_2(\bb Z_p) \times PGL_2(\bb Z_l)$.
\end{definition}

Since the action of $\Gamma$ on $\Delta$ is vertex transitive, every 
maximal abelian subgroup $\Gamma_0 \cong \bb Z^2$ is conjugate in $\Gamma$ to a period subgroup.
We want to show that $n(x)$ determines when $Z_{\Gamma}(\psi(x))$ is a period subgroup of $\Gamma$.

Recall that $\Gamma$ is generated by 
free groups $\Gamma_p$, $\Gamma_l$, of ranks $(p+1)/2$, $(l+1)/2$ respectively.
If $\gamma \in \Gamma$, let $\ell(\gamma)$ denote the natural word length of $\gamma$, in terms of the
generators of $\Gamma_p$, $\Gamma_l$. The condition $\ell(\gamma^2)=2\ell(\gamma)$, which is used in the
next lemma, is equivalent to the assertion that $\gamma$ has an axis containing $O$, upon which $\gamma$
acts by translation.

\begin{lemma}\label{hv} Let $a = \psi(x) \in \Gamma_p-\{1\}$ and let $n=n(x)$. 
The following statements are equivalent.
\begin{itemize}
\item[(a)] $p \nmid n$;
\item[(b)] $\ell(a^2)=2\ell(a)$;
\item[(c)] $\left(\frac{-n}{p}\right)=1$.
\end{itemize}
Similar equivalent assertions hold, if $p$ is replaced by $l$.
\end{lemma}

\noindent 
Before giving the proof, we note that
$$
\left(\frac{-n}{p}\right)= 
\begin{cases}
\left(\frac{n}{p}\right), & \text{if } p \equiv 1 \pmod 4 \,, \\
-\left(\frac{n}{p}\right), & \text{if } p \equiv 3 \pmod 4  \,.
\end{cases} 
$$

\begin{proof} $(a)\Leftrightarrow (b)$. The idea for this comes from the proof of \cite[Proposition 3.15]{moz}.
Write $x$ as in (\ref{x})
with $|x|^2 = x_0^2 + n z_0^2 = p^r$, $r > 0$.
Extracting a common factor, if necessary, we may assume $\gcd(x_0 , z_0) = 1$.
This means that $r = \ell(a)$ \cite[Corollary~3.11(4), Theorem~3.30(1)]{rat}.

Suppose that $p \nmid n$.
To prove $\ell(a^2)=2\ell(a)$ we must show that $p$ does not divide $x^2$. 
Now if $p$ divides 
$$x^2 = (x_0^2 - n z_0^2) + 2 x_0 z_0(c_1i+c_2j+c_3k)\,,$$
then $p$ divides the real part $x_0^2 - nz_0^2$.
Therefore $p$ divides $x_0$ (since $p$ divides $p^r = x_0^2+nz_0^2$).
But this implies that  $p$ divides $z_0$, since $p \nmid n$. This contradicts
the assumption $\gcd(x_0 , z_0) = 1$.

Conversely, suppose that $\ell(a^2)=2\ell(a)$. 
If $p$ divides $n$, then $p$ divides $x_0$ (since $p$ divides $x_0^2+nz_0^2$).
Therefore $p$ divides the real and imaginary parts of
$x^2=(x_0^2-nz_0^2)+2x_0z_0(c_1i+c_2j+c_3k)$. But this implies that 
$\ell(a^2)<2r$, a contradiction. 

$(a)\Leftrightarrow (c)$.
Suppose that $p \nmid n$. 
Note that $p$ does not divide $z_0$: otherwise $p$ also divides $x_0$.
It follows that $z_0$ has a multiplicative inverse$\pmod p$.
That is, one can choose $t \in \bb Z$ such that $z_0t\equiv 1 \pmod p$. Then
$$0\equiv (x_0^2+nz_0^2)t^2\equiv x_0^2t^2+n \pmod p \,.$$
Since $p \nmid n$, this means that $\left(\frac{-n}{p}\right)=1$.
The converse is obvious.
\end{proof}

\begin{lemma}\label{necessary}
If $\Gamma_0 \cong \bb Z^2$ is a period subgroup of $\Gamma$ and $n=n(\Gamma_0)$, then 
$\left(\frac{-n}{p}\right)=\left(\frac{-n}{l}\right)=1$.
\end{lemma}

\begin{proof}
The group $\Gamma_0$ acts cocompactly by translation on the apartment
$\cl A_{\Gamma_0}$ containing the vertex $O$.
It follows that $\Gamma_0$ contains elements $a\in\Gamma_p-\{1\}$, $b\in\Gamma_l-\{1\}$.
These elements act freely by translation on the apartment, and so
$\ell(a^2)=2\ell(a)$, $\ell(b^2)=2\ell(b)$. Therefore 
$\left(\frac{-n}{p}\right)=\left(\frac{-n}{l}\right)=1$, by Lemma \ref{hv}. 
\end{proof}

\begin{lemma}\label{moz}
If $\gamma = \psi(x) \in \Gamma - (\Gamma_p \cup \Gamma_l)$ and $\gcd(n(x), pl)=1$, then $Z_{\Gamma}(\gamma)$
is a period subgroup of $\Gamma$.
\end{lemma}

\begin{proof}
Let $x=x_0+z_0(c_1i+c_2j+c_3k)$ as in (\ref{x})
and
$n=n(x)= c_1^2+c_2^2+c_3^2$.
We may assume $\gcd(x_0 , z_0) = 1$ and
$|x|^2=x_0^2+nz_0^2 = p^rl^s$, where $r,s\ge 1$ because $\psi(x) \not\in\Gamma_p\cup\Gamma_l$.

The assumption $\gcd(n, pl)=1$ implies  
that $\gcd(x_0z_0, pl) = 1$. For example, if  $p \mid x_0$ then $p \mid z_0$,
since $p \mid (x_0^2+nz_0^2)$ and $p \nmid n$. This contradicts $\gcd(x_0 , z_0) = 1$. Similarly
$p \nmid z_0$.
It follows from the ``if'' part of the proof of \cite[Proposition 3.15]{moz} 
(and an obvious generalization to the cases where $p \equiv 3 \pmod 4$ or $l \equiv 3 \pmod 4$)
that $\gamma = \psi(x)$ lies in an abelian subgroup  $\Gamma_0$ of $\Gamma$, with $\Gamma_0 \cong \bb Z^2$.
The same proof also shows that $\Gamma_0$ acts cocompactly by translation on an apartment
$\cl A$ containing $O$.
(The essential point in the proof of Mozes is that $\ell(\gamma^2)=2\ell(\gamma)$.) 
However, $Z_{\Gamma}(\gamma)$ is the unique maximal abelian subgroup
containing $\Gamma_0$.
Therefore $Z_{\Gamma}(\gamma)$ acts cocompactly by translation on the apartment $\cl A$,
by the uniqueness assertion in \cite[6.8]{pra}.
In other words, $Z_{\Gamma}(\gamma)$ is a period subgroup of $\Gamma$.
\end{proof}

Now we can describe the period subgroups of $\Gamma$.

\begin{proposition}\label{per}
Let $\Gamma_0$ be a maximal abelian subgroup of $\Gamma$, and let $n=n(\Gamma_0)$.
Then $\Gamma_0$ is a period subgroup if and only if 
$\left(\frac{-n}{p}\right)=\left(\frac{-n}{l}\right)=1$.
\end{proposition}

Before proceeding with the proof, we introduce some notation. There is a
canonical Cartan subgroup $C$ of $G = PGL_2(\bb Q_p) \times PGL_2(\bb Q_l)$ defined by
$$C=\left(
\begin{pmatrix}
* & 0 \\
0  & * \\
\end{pmatrix},
\begin{pmatrix}
* & 0 \\
0  & * \\
\end{pmatrix}
\right)
\cap G\,.
$$
The group $C$ acts by translation on an apartment $\cl A$, which contains 
the vertex $O$ whose stabilizer in $G$ is $PGL_2(\bb Z_p) \times PGL_2(\bb Z_l)$.
The action of $C$ is transitive on the vertices of $\cl A$.

\begin{proof}[Proof of Proposition~\ref{per}]
In view of Lemma \ref{necessary}, it suffices to show that
$\left(\frac{-n}{p}\right)=\left(\frac{-n}{l}\right)=1$
implies that $\Gamma_0$ is a period subgroup.
Suppose therefore that $\left(\frac{-n}{p}\right)=\left(\frac{-n}{l}\right)=1$. Then $\gcd(n, pl)=1$. 
The result will therefore follow from Lemma \ref{moz}, if we can show
that $\Gamma_0$ is not contained in $\Gamma_p\cup\Gamma_l$.
By symmetry it is enough to prove that if $\Gamma_0$ contains an element
$b = \psi(y) \in \Gamma_l-\{1\}$, then it also contains an element $a = \psi(x) \in  \Gamma_p-\{1\}$.
For then the element $ba$ does not lie in $\Gamma_p\cup\Gamma_l$.

Write $y=y_0 + z_0(c_1 i + c_2 j + c_3 k)$,
where $c_1, c_2, c_3\in\bb Z$ are relatively prime
and $n=n(y)= c_1^2+c_2^2+c_3^2$.
The quaternion $y$ represents the element $b$ of $\Gamma_l$ of word length $\ell(b)=s>0$.
By Lemma \ref{hv}, $b$ acts by translation of distance $s$ along an axis $L_b$ containing $O$.

The element of $GL_2(\bb Q_p) \times GL_2(\bb Q_l)$ corresponding to $y$ in the formula 
(\ref{psimap}) has eigenvalues $y_0 \pm z_0\sqrt{-n}$. 
The assumption $\left(\frac{-n}{p}\right)=\left(\frac{-n}{l}\right)=1$
implies that $\sqrt{-n}$ exists in both
$\bb Q_p$ and $\bb Q_l$ and therefore that $b$ is diagonalizable in $G$. 
In other words, there exists an element $h \in G$ such that $h^{-1}bh \in C$.

The group $hCh^{-1}$ acts by translation on the apartment $h\cl A$. 
Also the element $b\in hCh^{-1}\cap \Gamma_l$ acts by translation on the apartment $h\cl A$, in a direction which will be
called ``vertical''. 
Now $h\cl A$ necessarily contains the axis $L_b$ of $b$, by \cite[Theorem II.6.8 (3)]{bh}.
In particular, $O\in h\cl A$. 
 
Choose $g\in hCh^{-1}$ to act on $h\cl A$ by horizontal translation. 
Consider the horizontal strip $H$ in $h\cl A$ obtained by translating the vertical segment $[O, bO]$. 
\refstepcounter{picture}
\begin{figure}[htbp]
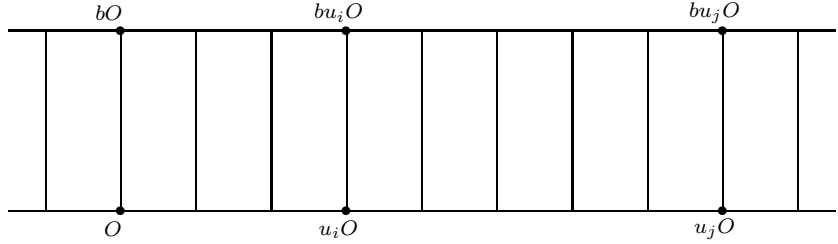
\label{strip}
\hfil
\centerline{
\beginpicture
\setcoordinatesystem units <0.5cm, 1.2cm>  
\setplotarea  x from -10 to 10,  y from -1 to 1.2
\put {$_{O}$}     at    -8.2 -1.2
\put {$_{bO}$}     at    -8.3 1.2
\put {$_{u_iO}$}     at   -2.2 -1.2
\put {$_{bu_iO}$}    at   -2.2  1.2
\put {$_{u_jO}$}     at   7.8 -1.2
\put {$_{bu_jO}$}    at   7.8  1.2
\put {$_{\bullet}$} at -8 -1
\put {$_{\bullet}$} at -8  1
\put {$_{\bullet}$} at -2 -1
\put {$_{\bullet}$} at -2  1
\put {$_{\bullet}$} at  8 -1
\put {$_{\bullet}$} at   8  1
\putrule from -11 -1 to 11 -1
\putrule from -11  1 to 11  1
\putrule from 0 -1 to 0 1
\putrule from  2  -1 to 2  1
\putrule from  4  -1 to 4  1
\putrule from  6  -1 to 6  1
\putrule from  8  -1 to 8  1
\putrule from  10 -1 to 10  1
\putrule from  -2  -1 to -2  1
\putrule from  -4  -1 to -4  1
\putrule from  -6  -1 to -6  1
\putrule from  -8  -1 to -8  1
\putrule from  -10 -1 to -10  1
\endpicture
}
\hfil
\caption{The horizontal strip $H$.}
\end{figure}

Since $\Gamma$ acts freely and transitively on the vertices of $\Delta$, each vertical segment
$g^i[O,bO]$ of $H$ lies in the $\Gamma$-orbit of precisely one segment of the form $[O, \gamma O]$,
$\gamma \in \Gamma_l$, $\ell(\gamma)=s$. Moreover, there are only finitely many such segments
$[O, \gamma O]$.

If $i>0$ then $g^iO = u_iO$, for some $u_i\in\Gamma_p - \{ 1 \}$.
Since $b$ and $g$ commute, we have $g^ibO = bg^iO = bu_iO$.
That is, $g^i[O, bO]= [u_iO, bu_iO]$, which lies in the $\Gamma$-orbit
of the segment $[O, u_i^{-1}bu_iO]$.  
By the finiteness assertion in the preceding paragraph, there exist 
integers $j>i>0$
such that 
$$[O, u_i^{-1}bu_iO]=[O, u_j^{-1}bu_jO]\,.$$
By freeness of the action of $\Gamma$, 
$$u_i^{-1}bu_i = u_j^{-1}bu_j\,,$$
and $u_i\not=u_j$.
Therefore $ab=ba$, where $a=u_iu_j^{-1}\in  \Gamma_p-\{1\}$.
\end{proof}

A maximal abelian subgroup $\Gamma_0$ of $\Gamma$ may be isomorphic to $\bb Z$. 
Here is a way of providing some examples. 

\begin{corollary}\label{rank1} Suppose that $a\in \Gamma_p-\{1\}$, and $n=n(a)$ satisfies
$$\left(\frac{-n}{p}\right) = 1, \quad \left(\frac{-n}{l}\right) = -1\,.$$
Then $Z_{\Gamma}(a)<\Gamma_p$ is a maximal abelian subgroup of $\Gamma$, and
$Z_{\Gamma}(a)\cong \bb Z$.\\
A similar assertion applies to elements of $\Gamma_l-\{1\}$.
\end{corollary}

\begin{proof}The hypothesis implies that $\gcd(n, pl)=1$. 
If $Z_{\Gamma}(a)\not\subset\Gamma_p$, then $Z_{\Gamma}(a)$ contains an element 
$\gamma \not\in\Gamma_p\cup \Gamma_l$. Therefore $Z_{\Gamma}(a) = Z_{\Gamma}(\gamma)$ is a period group,
by Lemma \ref{moz}.
But this implies $\left(\frac{-n}{l}\right) = 1$, by Proposition \ref{per},
-- a contradiction.
\end{proof}

\begin{figure}[htbp]
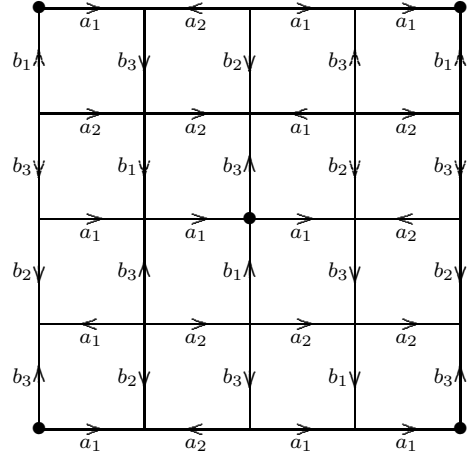
\label{periodic35}
\hfil
\centerline{
\beginpicture
\setcoordinatesystem units <0.7cm, 0.7cm>  
\setplotarea  x from 0 to 10,  y from 0 to 8
\putrule from 0 0 to 8 0
\putrule from 0 2 to 8 2
\putrule from 0 4 to 8 4
\putrule from 0 6 to 8 6
\putrule from 0 8 to 8 8
\arrow <6pt> [.3,.67] from    0.9  0 to   1.2 0
\arrow <6pt> [.3,.67] from    3.1  0 to   2.8 0
\arrow <6pt> [.3,.67] from    4.9  0 to   5.2 0
\arrow <6pt> [.3,.67] from    6.9  0 to   7.2 0
\put {$_{a_1}$}      at     1 -0.3
\put {$_{a_2}$}      at     3 -0.3
\put {$_{a_1}$}      at     5 -0.3
\put {$_{a_1}$}      at     7 -0.3
\arrow <6pt> [.3,.67] from    1.1  2 to   0.8 2
\arrow <6pt> [.3,.67] from    2.9  2 to   3.2 2
\arrow <6pt> [.3,.67] from    4.9  2 to   5.2 2
\arrow <6pt> [.3,.67] from    6.9  2 to   7.2 2
\put {$_{a_1}$}      at     1  1.7
\put {$_{a_2}$}      at     3  1.7
\put {$_{a_2}$}      at     5  1.7
\put {$_{a_2}$}      at     7  1.7
\arrow <6pt> [.3,.67] from    0.9  4 to   1.2 4
\arrow <6pt> [.3,.67] from    2.9  4 to   3.2 4
\arrow <6pt> [.3,.67] from    4.9  4 to   5.2 4
\arrow <6pt> [.3,.67] from    7.1  4 to   6.8 4
\put {$_{a_1}$}      at     1  3.7
\put {$_{a_1}$}      at     3  3.7
\put {$_{a_1}$}      at     5  3.7
\put {$_{a_2}$}      at     7  3.7
\arrow <6pt> [.3,.67] from    0.9  6 to   1.2 6
\arrow <6pt> [.3,.67] from    2.9  6 to   3.2 6
\arrow <6pt> [.3,.67] from    5.1  6 to   4.8 6
\arrow <6pt> [.3,.67] from    6.9  6 to   7.2 6
\put {$_{a_2}$}      at     1  5.7
\put {$_{a_2}$}      at     3  5.7
\put {$_{a_1}$}      at     5  5.7
\put {$_{a_2}$}      at     7  5.7
\arrow <6pt> [.3,.67] from    0.9  8 to   1.2 8
\arrow <6pt> [.3,.67] from    3.1  8 to   2.8 8
\arrow <6pt> [.3,.67] from    4.9  8 to   5.2 8
\arrow <6pt> [.3,.67] from    6.9  8 to   7.2 8
\put {$_{a_1}$}      at     1  7.7
\put {$_{a_2}$}      at     3  7.7
\put {$_{a_1}$}      at     5  7.7
\put {$_{a_1}$}      at     7  7.7
\putrule from 0 0 to 0 8
\putrule from 2 0 to 2 8
\putrule from 4 0 to 4 8
\putrule from 6 0 to 6 8
\putrule from 8 0 to 8 8
\arrow <6pt> [.3,.67] from   0  0.9   to  0 1.2 
\arrow <6pt> [.3,.67] from   0  3.1   to  0 2.8 
\arrow <6pt> [.3,.67] from   0  5.1   to  0 4.8
\arrow <6pt> [.3,.67] from   0  6.9   to  0 7.2 
\put {$_{b_3}$}      at   -0.3  1
\put {$_{b_2}$}      at   -0.3  3
\put {$_{b_3}$}      at   -0.3  5
\put {$_{b_1}$}      at   -0.3  7
\arrow <6pt> [.3,.67] from   2  1.1   to  2 0.8 
\arrow <6pt> [.3,.67] from   2  2.9   to  2 3.2 
\arrow <6pt> [.3,.67] from   2  5.1   to  2 4.8
\arrow <6pt> [.3,.67] from   2  7.1   to  2 6.8 
\put {$_{b_2}$}      at   1.7  1
\put {$_{b_3}$}      at   1.7  3
\put {$_{b_1}$}      at   1.7  5
\put {$_{b_3}$}      at   1.7  7
\arrow <6pt> [.3,.67] from   4  1.1   to  4 0.8 
\arrow <6pt> [.3,.67] from   4  2.9   to  4 3.2 
\arrow <6pt> [.3,.67] from   4  4.9   to  4 5.2
\arrow <6pt> [.3,.67] from   4  7.1   to  4 6.8 
\put {$_{b_3}$}      at   3.7  1
\put {$_{b_1}$}      at   3.7  3
\put {$_{b_3}$}      at   3.7  5
\put {$_{b_2}$}      at   3.7  7
\arrow <6pt> [.3,.67] from   6  1.1   to  6 0.8 
\arrow <6pt> [.3,.67] from   6  3.1   to  6 2.8 
\arrow <6pt> [.3,.67] from   6  5.1   to  6 4.8
\arrow <6pt> [.3,.67] from   6  6.9   to  6 7.2 
\put {$_{b_1}$}      at   5.7  1
\put {$_{b_3}$}      at   5.7  3
\put {$_{b_2}$}      at   5.7  5
\put {$_{b_3}$}      at   5.7  7
\arrow <6pt> [.3,.67] from   8  0.9   to  8 1.2 
\arrow <6pt> [.3,.67] from   8  3.1   to  8 2.8 
\arrow <6pt> [.3,.67] from   8  5.1   to  8 4.8
\arrow <6pt> [.3,.67] from   8  6.9   to  8 7.2 
\put {$_{b_3}$}      at   7.7  1
\put {$_{b_2}$}      at   7.7  3
\put {$_{b_3}$}      at   7.7  5
\put {$_{b_1}$}      at   7.7  7
\put {$\bullet$} at 0 0 
\put {$\bullet$} at 8 0 
\put {$\bullet$} at 4 4
\put {$\bullet$} at 0 8 
\put {$\bullet$} at 8 8 
\endpicture
}
\hfil
\caption{Part of a periodic apartment for $\Gamma_0 < \Gamma_{3,5}$\,.}
\end{figure}

\begin{example} \label{ex35}
Let $\Gamma = \Gamma_{3,5}$. This group has a presentation with generators 
$\{a_1, a_2, b_1, b_2, b_3\}$ and relators 
$$\{a_1 b_1 a_2 b_2, \, a_1 b_2 a_2 b_1^{-1}, \,
    a_1 b_3 a_2^{-1} b_1, \, a_1 b_3^{-1} a_1 b_2^{-1}, \, 
    a_1 b_1^{-1} a_2^{-1} b_3, \, a_2 b_3 a_2 b_2^{-1}\}\,,
$$
where
\begin{align}
a_1 &= \psi( 1 + j + k), & a_1^{-1} &= \psi( 1 - j - k), \notag \\    
a_2 &= \psi( 1 + j - k), & a_2^{-1} &= \psi( 1 - j + k), \notag \\
b_1 &= \psi( 1 + 2i), & b_1^{-1} &= \psi( 1 - 2i), \notag \\    
b_2 &= \psi( 1 + 2j), & b_2^{-1} &= \psi( 1 - 2j), \notag \\  
b_3 &= \psi( 1 + 2k), & b_3^{-1} &= \psi( 1 - 2k).\notag
\end{align}

The subgroup $\langle a_1 \rangle = Z_{\Gamma}(a_1) < \Gamma_3$ is maximal abelian in $\Gamma$
by Corollary~\ref{rank1},
since $n(a_1) = 2$, $\left(\frac{-2}{3}\right)=1$ and $\left(\frac{-2}{5}\right)=-1$.

The subgroup $\langle a_1 a_2^{-1}  a_1^2 \rangle = \langle \psi(-5-6i-2j+4k) \rangle$ is not maximal abelian. 
It is contained in the period subgroup
$$\Gamma_0=\langle a_1 a_2^{-1} a_1^2, \, b_3 b_2^{-1} b_3^{-1} b_1 \rangle \cong \bb Z^2\,.$$
Indeed, $n(\Gamma_0)= n(a_1 a_2^{-1} a_1^2) =14$, $\left(\frac{-14}{3}\right)=1$, $\left(\frac{-14}{5}\right)=1$.
Note that $b_3 b_2^{-1} b_3^{-1} b_1 = \psi(-11 + 18i + 6j - 12k)$. 
Part of the period lattice for $\Gamma_0$ is illustrated in Figure~\ref{periodic35}.
\end{example}

\refstepcounter{picture}

\begin{figure}[htbp]
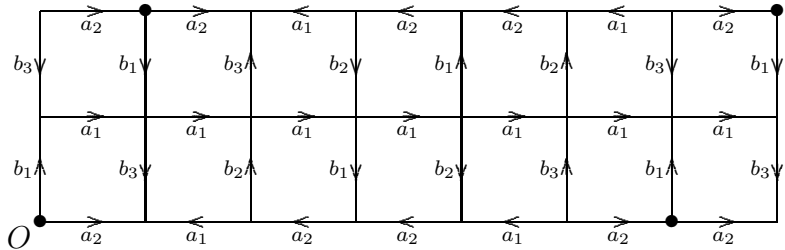
\label{periodicstrip35}
\hfil
\centerline{
\beginpicture
\setcoordinatesystem units <0.7cm, 0.7cm>  
\setplotarea  x from 0 to 14,  y from 0 to 4
\putrule from 0 0 to 14 0
\putrule from 0 2 to 14 2
\putrule from 0 4 to 14 4
\putrule from 0 0 to 0 4
\putrule from 2 0 to 2 4
\putrule from 4 0 to 4 4
\putrule from 6 0 to 6 4
\putrule from 8 0 to 8 4
\putrule from 10 0 to 10 4
\putrule from 12 0 to 12 4
\putrule from 14 0 to 14 4
\arrow <6pt> [.3,.67] from    0.9  0 to   1.2 0
\arrow <6pt> [.3,.67] from    3.1  0 to   2.8 0
\arrow <6pt> [.3,.67] from    5.1  0 to   4.8 0
\arrow <6pt> [.3,.67] from    7.1  0 to   6.8 0
\arrow <6pt> [.3,.67] from    9.1  0 to   8.8 0
\arrow <6pt> [.3,.67] from   10.9  0 to  11.2 0
\arrow <6pt> [.3,.67] from   12.9  0 to  13.2 0
\put {$_{a_2}$}  at   1  -0.3
\put {$_{a_1}$}  at   3  -0.3
\put {$_{a_2}$}  at   5  -0.3
\put {$_{a_2}$}  at   7  -0.3
\put {$_{a_1}$}  at   9  -0.3
\put {$_{a_2}$}  at  11  -0.3
\put {$_{a_2}$}  at  13  -0.3
\arrow <6pt> [.3,.67] from    0.9  2 to   1.2 2
\arrow <6pt> [.3,.67] from    2.9  2 to   3.2 2
\arrow <6pt> [.3,.67] from    4.9  2 to   5.2 2
\arrow <6pt> [.3,.67] from    6.9  2 to   7.2 2
\arrow <6pt> [.3,.67] from    8.9  2 to   9.2 2
\arrow <6pt> [.3,.67] from   10.9  2 to  11.2 2
\arrow <6pt> [.3,.67] from   12.9  2 to  13.2 2
\put {$_{a_1}$}  at   1  1.7
\put {$_{a_1}$}  at   3  1.7
\put {$_{a_1}$}  at   5  1.7
\put {$_{a_1}$}  at   7  1.7
\put {$_{a_1}$}  at   9  1.7
\put {$_{a_1}$}  at  11  1.7
\put {$_{a_1}$}  at  13  1.7
\arrow <6pt> [.3,.67] from    0.9  4 to   1.2 4
\arrow <6pt> [.3,.67] from    2.9  4 to   3.2 4
\arrow <6pt> [.3,.67] from    5.1  4 to   4.8 4
\arrow <6pt> [.3,.67] from    7.1  4 to   6.8 4
\arrow <6pt> [.3,.67] from    9.1  4 to   8.8 4
\arrow <6pt> [.3,.67] from   11.1  4 to  10.8 4
\arrow <6pt> [.3,.67] from   12.9  4 to  13.2 4
\put {$_{a_2}$}  at   1  3.7
\put {$_{a_2}$}  at   3  3.7
\put {$_{a_1}$}  at   5  3.7
\put {$_{a_2}$}  at   7  3.7
\put {$_{a_2}$}  at   9  3.7
\put {$_{a_1}$}  at  11  3.7
\put {$_{a_2}$}  at  13  3.7
\arrow <6pt> [.3,.67] from   0  0.9   to  0  1.2 
\arrow <6pt> [.3,.67] from   0  3.1   to  0  2.8 
\put {$_{b_1}$}  at    -0.3  1
\put {$_{b_3}$}  at    -0.3  3
\arrow <6pt> [.3,.67] from   2  1.1   to  2  0.8 
\arrow <6pt> [.3,.67] from   2  3.1   to  2  2.8 
\put {$_{b_3}$}  at    1.7  1
\put {$_{b_1}$}  at    1.7  3
\arrow <6pt> [.3,.67] from   4  0.9   to  4  1.2 
\arrow <6pt> [.3,.67] from   4  2.9   to  4  3.2 
\put {$_{b_2}$}  at    3.7  1
\put {$_{b_3}$}  at    3.7  3
\arrow <6pt> [.3,.67] from   6  1.1   to  6  0.8 
\arrow <6pt> [.3,.67] from   6  3.1   to  6  2.8 
\put {$_{b_1}$}  at    5.7  1
\put {$_{b_2}$}  at    5.7  3
\arrow <6pt> [.3,.67] from   8  1.1   to  8  0.8 
\arrow <6pt> [.3,.67] from   8  2.9   to  8  3.2 
\put {$_{b_2}$}  at    7.7  1
\put {$_{b_1}$}  at    7.7  3
\arrow <6pt> [.3,.67] from   10  0.9   to  10  1.2 
\arrow <6pt> [.3,.67] from   10  2.9   to  10  3.2 
\put {$_{b_3}$}  at    9.7  1
\put {$_{b_2}$}  at    9.7  3
\arrow <6pt> [.3,.67] from   12  0.9   to  12  1.2 
\arrow <6pt> [.3,.67] from   12  3.1   to  12  2.8 
\put {$_{b_1}$}  at    11.7  1
\put {$_{b_3}$}  at    11.7  3
\arrow <6pt> [.3,.67] from   14  1.1   to  14  0.8 
\arrow <6pt> [.3,.67] from   14  3.1   to  14  2.8 
\put {$_{b_3}$}  at   13.7  1
\put {$_{b_1}$}  at   13.7  3
\put {$\bullet$} at 0 0
\put {$O$}       at -0.4 -0.3 
\put {$\bullet$} at 12 0
\put {$\bullet$} at 2 4
\put {$\bullet$} at 14 4
\endpicture
}
\hfil
\caption{Part of a periodic horizontal strip\,.}
\end{figure}

\begin{example}
Let $\Gamma=\Gamma_{3,5}$.
Consider $b_1 a_1 b_1^{-1} = \psi(5-7j+k)$.
By Example \ref{ex35}, $\langle a_1 \rangle$ is maximal abelian in $\Gamma$.
Therefore so also is $\Gamma_0=\langle  b_1 a_1 b_1^{-1} \rangle = b_1 \langle a_1 \rangle b_1^{-1}$.
Now $\gamma = b_1 a_1^6 b_1^{-1} = a_2 a_1^{-1} a_2^{-2} a_1^{-1} a_2 = \psi(5(23+14j-2k)) = \psi(x) \in \Gamma_{3}$,
with $|x|^2=5^2.3^6$.  Also $n(x)=n(\Gamma_0) = 50$, $\left(\frac{-50}{3}\right)=1$ and $\left(\frac{-50}{5}\right)=0$. 
There is a periodic horizontal strip of height 2 (Figure \ref{periodicstrip35}), 
upon which $\gamma$ acts by translation. This strip is the union of the axes of $\gamma$.
\end{example}

\refstepcounter{picture}

\begin{example}
Let $\Gamma = \Gamma_{3,5}$.
Conjugating the period subgroup 
$\langle a_1 a_2^{-1} a_1^2, \, b_3 b_2^{-1} b_3^{-1} b_1 \rangle$
of Example~\ref{ex35} by $a_2$ gives the group
\begin{align}
\Gamma_0 &= \langle a_2 a_1 a_2^{-1} a_1^2 a_2^{-1}, \, a_2 b_3 b_2^{-1} b_3^{-1} b_1 a_2^{-1} \rangle =
\langle a_2 a_1 a_2^{-1} a_1^2 a_2^{-1}, \, b_2 b_1^{-1} b_2^2 \rangle \notag \\
&= \langle \psi(-15+10i+2j+20k), \, \psi(-11-10i-2j-20k) \rangle \cong \bb Z^2\,, \notag
\end{align}
which is not a period subgroup since $n(\Gamma_0)=126$,
$\left(\frac{-126}{5}\right)=1$ and $\left(\frac{-126}{3}\right)=0$.
\end{example}

One could conjecture that every maximal abelian subgroup of $\Gamma$ is conjugate to either a
period subgroup or to a subgroup of $\Gamma_p$ or $\Gamma_l$.
The next example shows that this conjecture is not true.
We need the following definition and Lemma~\ref{mlemma}:

If $x = x_0 + x_1 i + x_2 j + x_3 k \in \bb H(\bb Z)$, 
let $m(x) = |x|^2 - \Re(x)^2 = x_1^2 + x_2^2 + x_3^2$,
where $\Re(x) = x_0$ denotes the real part of $x$.
Observe that $m(x) = \lambda^2 n(x)$ for some integer $\lambda$.

\begin{lemma}\label{mlemma}
Let $x, y \in \bb H(\bb Z)$, then $m(xy\ovl x) = (|x|^2)^2 m(y)$.
\end{lemma}

\begin{proof}
Using the rules $\Re(xy) = \Re(yx)$ and $|xy|^2 = |x|^2 |y|^2$, we conclude
$$m(xy\overline{x}) = |xy\ovl x|^2 - \Re(xy\ovl x)^2 = (|x|^2)^2 |y|^2 - (|x|^2 \Re(y))^2  = (|x|^2)^2 m(y)\, .$$
\end{proof}

\begin{figure}[htbp]
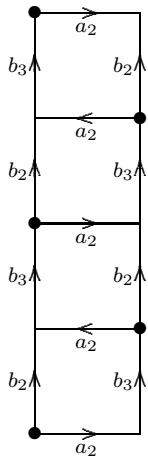
\label{periodica2b3}
\hfil
\centerline{
\beginpicture
\setcoordinatesystem units <0.7cm, 0.7cm>  
\setplotarea  x from 0 to 2,  y from 0 to 8
\putrule from 0 0 to 0 8
\putrule from 2 0 to 2 8
\putrule from 0 0 to 2 0
\putrule from 0 2 to 2 2
\putrule from 0 4 to 2 4
\putrule from 0 6 to 2 6
\putrule from 0 8 to 2 8
\arrow <6pt> [.3,.67] from    0.9  0 to   1.2 0
\put {$_{a_2}$}      at   1 -0.3
\arrow <6pt> [.3,.67] from    1.1  2 to   0.8 2
\put {$_{a_2}$}      at   1  1.7
\arrow <6pt> [.3,.67] from    0.9  4 to   1.2 4
\put {$_{a_2}$}      at   1  3.7
\arrow <6pt> [.3,.67] from    1.1  6 to   0.8 6
\put {$_{a_2}$}      at   1  5.7
\arrow <6pt> [.3,.67] from    0.9  8 to   1.2 8
\put {$_{a_2}$}      at   1  7.7
\arrow <6pt> [.3,.67] from   0  0.9   to  0 1.2 
\arrow <6pt> [.3,.67] from   0  2.9   to  0 3.2 
\arrow <6pt> [.3,.67] from   0  4.9   to  0 5.2 
\arrow <6pt> [.3,.67] from   0  6.9   to  0 7.2
\put {$_{b_2}$}      at    -0.3 1 
\put {$_{b_3}$}      at    -0.3 3 
\put {$_{b_2}$}      at    -0.3 5 
\put {$_{b_3}$}      at    -0.3 7 
\arrow <6pt> [.3,.67] from   2  0.9   to  2 1.2 
\arrow <6pt> [.3,.67] from   2  2.9   to  2 3.2 
\arrow <6pt> [.3,.67] from   2  4.9   to  2 5.2 
\arrow <6pt> [.3,.67] from   2  6.9   to  2 7.2
\put {$_{b_3}$}      at    1.7 1 
\put {$_{b_2}$}      at    1.7 3 
\put {$_{b_3}$}      at    1.7 5 
\put {$_{b_2}$}      at    1.7 7 
\put {$\bullet$} at 0 0 
\put {$\bullet$} at 2 2 
\put {$\bullet$} at 0 4
\put {$\bullet$} at 2 6 
\put {$\bullet$} at 0 8 
\endpicture
}
\hfil
\caption{Part of a periodic vertical strip\,.}
\end{figure}

\begin{example}
Let $\Gamma = \Gamma_{3,5}$ and $a_2 b_3 = \psi(3+2i+j+k)$.
The group $\Gamma_0 = Z_{\Gamma}(a_2 b_3)$ is a maximal abelian subgroup of $\Gamma$ such that $n(\Gamma_0) = 6$.
We fix any element $\gamma = \psi(x) \in \Gamma$.

The maximal abelian subgroup $\gamma \Gamma_0 \gamma^{-1}$ is not a subgroup of $\Gamma_3$ or $\Gamma_5$,
since $\gamma a_2 b_3 \gamma^{-1} \in \gamma \Gamma_0 \gamma^{-1}$ is the $\psi$-image of $x(3+2i+j+k)\ovl x$
whose norm is a product of an odd power of $3$ and an odd power of $5$.

We claim that $\gamma \Gamma_0 \gamma^{-1}$ is not a period subgroup. If $|x|^2 = 3^r 5^s$, $r,s \geq 0$, then
by Lemma~\ref{mlemma}
$$(3^r 5^s)^2.6 = m(x(3+2i+j+k)\ovl x) = \lambda^2 n(\gamma \Gamma_0 \gamma^{-1})$$
for some integer $\lambda$. It follows that $3 \mid n(\gamma \Gamma_0 \gamma^{-1})\,$,
in particular 
$$\left(\frac{-n(\gamma \Gamma_0 \gamma^{-1})}{3}\right)=0$$ 
and Proposition~\ref{per} proves the claim.

Since any maximal abelian subgroup of rank 2 is conjugate to a period subgroup, it also follows
that $\Gamma_0 \cong \bb Z$. 
See Figure~\ref{periodica2b3} for a periodic vertical strip of width $1$ 
which is globally invariant under the action of $a_2 b_3$. Note that $(a_2 b_3)^2 = b_2 b_3$.
Therefore $a_2 b_3$ acts upon the strip by glide reflection and the unique axis of $a_2 b_3$ is
the vertical central line of the strip.
\end{example}

\refstepcounter{picture}

It is well-known that period subgroups in $\Gamma$ always exist.
See for example \cite[Proposition~4.2]{rat} for an elementary proof of this fact, using
doubly periodic tilings of the Euclidean plane by unit squares.
We mention a corollary of this in terms of integer quaternions.

\begin{corollary} \label{corper}
Given any pair $(p,l)$ of distinct odd primes,
there are $x,y \in \tilde \Gamma$ and $1 \leq r \leq 4(p+1)^2 (l+1)^2$ such that $xy = yx$ and 
$$|x|^2 = p^r, |y|^2 = l^r, \left( \frac{-n(x)}{p} \right) = \left( \frac{-n(y)}{l} \right) = 1\,.$$ 
\end{corollary}

The integer $r$ in this corollary comes from the 
constructive proof of \cite[Proposition~4.2]{rat}, 
and its upper bound is certainly not optimal.
In fact, if $p,l \equiv 1 \pmod 4$, 
there is a direct proof of Corollary~\ref{corper} (with $r = 1$),
applying the Two Square Theorem.

\end{document}